\documentclass[10pt,a4paper,fleqn]{article}
\usepackage[utf8]{inputenc}
\usepackage{amsmath}
\usepackage{amsfonts}
\usepackage{amssymb}
\usepackage[dvipsnames]{xcolor}
\usepackage{graphicx}
\usepackage{subcaption}
\usepackage{float}
\usepackage{enumitem}

\numberwithin{equation}{section}

\author{Giulio Morpurgo}
\title{A simple estimate of the number of Goldbach pairs for $N$}
\begin{document}
\maketitle
\begin{abstract}
Using the fact that the number of combinations $p_{1}$, $p_{2}$, where $p_{1}$ and $p_{2}$ are odd primes, with $p_{1} \leq p_{2}$ and $p_{1} + p_{2} \leq 2N$ is equal to the total number of Goldbach pairs for all the even integers from 6 to $2N$, we derive in a simple way an accurate estimate of the number of Goldbach pairs for $N$.
\end{abstract}

\section{Introduction}

We call ``Goldbach pair" for $N$ a pair of prime numbers $N-m$, $N+m$ (whose sum is $2N$). 
\\
The existence of a Goldbach pair $N-m$, $N+m$ is trivially equivalent to the existence of two primes at interval $2m$. Every pair of odd primes $p_1, p_2$ is a Goldbach pair for $(p_1 + p_2)/2$. If we exclude the Goldbach pairs $N, N$, where $N$ is a prime, the set of all remaining Goldbach pairs for all the $N$ values up to $N_{max}$ is identical to the set of all pairs of primes $(p_1$,\;$p_2)$ an any even interval from 2 to $2N_{max} - 6$ \footnote{as 3 is the lowest odd prime}, such that  $p_1 < p_2$ and their sum is not larger than $2N_{max}$. 
\\
Also, it is a known fact (see, for example, \cite{Hardy_Littlewood_1}, \cite{Fl_Ro}) that for two integers $N_{1}$ and $N_{2}$ of similar size (for instance $N_{1}=N$ and $N_{2}=N+1$) the ratio between the number of Goldbach pairs for $N_{1}$ and the number of Goldbach pairs for $N_{2}$ mainly depends on which odd primes appear in the factorizations of  $N_{1}$ and $N_{2}$ respectively.
\\
From these two facts it is possible to formulate a simple and accurate estimate of the number of Goldbach pairs for $N$.

\section{Counting the total number of Goldbach pairs}
In this section, we first introduce a discrete ``counting" function, defined on the integers, that counts the total number of Goldbach pairs for all the integers up to a value $M$. Then we search a way to approximate this discrete function using a continuous one. Finally we show how this continuous function can be used to get an estimate for the number of Goldbach pairs for any value of $N$.  
\subsection{\normalsize A discrete function for the total number of pairs of odd primes whose sum is not larger than $M$}
We can define a discrete function $G_{tot}(M)$. The value of $G_{tot}(M)$ is the total number of Goldbach pairs with sum $\leq M$. This function is defined over the integers, and is monotonically increasing ($G_{tot}(M-1) \leq G_{tot}(M)$. When $M$ increases by 1 (from $M-1$ to $M$), the value of $G_{tot}(M)$ increases by the number of odd prime pairs summing to $M$. $G_{tot}(M) - G_{tot}(M-1)$ can only be larger than zero if $M$ is even; if $M$ is odd there is no new Goldbach pair, and the function keeps its previous value. The first few values of $G_{tot}(M)$, for $M$ from 6 to 20, are 1, 1, 2, 2, 4, 4, 5, 5, 7, 7, 9, 9, 11, 11, 13.
\subsection{\normalsize A continuous approximation for $G_{tot}(M)$}
\label{sec-approx}
The number of Goldbach pairs whose sum is not larger than $2N$ depends on the number of primes not larger than $2N-3$. Let be $P$ the number of odd primes up to $2N-3$. Among these $P$ prime numbers, one can select pairs for which the sum is below $2N$, and pairs for which the sum is above $2N$. If we make the assumptions that {\bf a}) the number of primes below $N$ is not too different than the number of primes from $N$ to $2N$, and {\bf b}) that there is not too much bias in the overall way these primes are distributed, then a reasonable first approximation for the total number of pairs satisfying the condition is given by
\begin{equation}
\label{eq:approx}
G_{tot}(2N) \approx P^{2}/4
\end{equation}
We know that in general the two assumptions a) and b) are not correct. If we look at the interval $0..2N$, the number of primes in the first half of the interval is larger than the number in the second half, and the density is much higher at the beginning. But, when $2N$ gets very large, at least the relative error due to assumption a) decreases. Even better, we can choose to only look at an interval in which a) and b) are asymptotically satisfied: the ``reduced" range. 
\\
In general, we look for Goldbach pairs $N-m$, $N+m$ with $m$ going from 0 to $N-3$ (``full" range). Within the full range, assumptions a) and b) are not correct, and the density of primes ranges over orders of magnitude. But if instead of allowing $m$ to run all the way up to $N-3$ we force it to stay below $(\sqrt{2}-1)N$, the primes composing our pairs range from $0.5858N$ to $1.4142N$ (``reduced" range). Within this range, the density of prime numbers goes from $1/(\log{N}-0.523)$ to $1/(\log{N}+0.347)$. When $N$ is very large, the density of these primes can be considered as almost constant. Later we will compare our estimates with the real results for both the full and the reduced range; we will see how the estimate is much more precise when applied to the reduced range.
\\
For now, we continue like if a) and b) were correct. A quantity we need to progress from equation \ref{eq:approx} is $P$, the number of primes up to $2N$. We can get an approximate value for $P$ using the following bound, proved by Pierre Dusart \cite{Dusart_2}:
\begin{equation}
\pi(x) < \frac{x}{\log{x}} \cdot \left(1 + \frac{1}{\log{x}} + \frac{2.51}{(\log{x})^2}\right)  \verb+ for +x \geq 355991  \label{eq:Dusart}
\end{equation}
This not only gives a max. bound for P; it also returns a value whose relative error compared to the exact one, $\pi(x)$, is small.\footnote{in the range from 5 millions to 10 billions, the value exceeds $\pi(x)$ only by less than four parts over ten thousand; for higher $x$ the relative precision increases.}
\\
We now plug in Dusart's result into our approximation \ref{eq:approx}, to obtain
\begin{equation}
\label{eq:approx2}
G_{tot}(x) \approx \frac{1}{4} \left( \frac{x}{\log{x}} \cdot \left(1 + \frac{1}{\log{x}} + \frac{2.51}{(\log{x})^2}\right) \right)^{2}
\end{equation} 
If we replace $x$ by $2N$, we get the value we were looking for. In general, we can look at the term on the right as a formula that provides meaningful values when we replace x with an integer ($> 355991$). Or else, we could treat it as a continuous function over the reals (even if the values mean something only when $x$ is an integer).  The value of this continuous function at x is (approximatively) the sum of the number of Goldbach pairs for all the values up to x. Let us rename this function as {\bf $\boldsymbol{ g_{tot}(x)}$}. In other words, {\bf $\boldsymbol {g_{tot}(x)}$  approximates the integral of an unknown ``function" that counts the Goldbach pairs for an integer N} (Let's call this unknown function {\bf g(x)}). 
\\{\bf If $\boldsymbol{g_{tot}(x)}$ is related to the integral of g(x), its derivative should be related to g(x)}.
\\
This is what the derivative of $\boldsymbol {g_{tot}(x)}$ looks like\footnote{as simplified by ``derivative-calculator.net"}
\begin{equation}
\label{eq:derivative}
\boldsymbol{g_{tot}'}(x) = \frac{1}{4} \cdot \frac{x \cdot (100 \log^{2}(x) + 100 \log(x) +251) \cdot (100 \log^{3}(x) + 51 \log(x) -753)}{5000 \cdot \log^{7}(x)}
\end{equation}
In what follows, we investigate the relation between $g_{tot}'(x)$ and $g(x)$, and show how an accurate approximation of $g(x)$ can be obtained from this relation. 

\section{Number of Goldbach pairs and factorization of N}
Before being able to use equation \ref{eq:derivative} and estimate the number of Goldbach pairs for a given N, there is one important point we have to address. We mentioned in the introduction that the number of Goldbach pairs {\em (N-m, N+m)} for a given N strongly depends on the factorization of $N$. But this fact does not show up in the above analysis. In fact, while the real function $G_{tot}(N)$ increases in steps whose values correspond to the number of Goldbach pairs for $N$, the approximate function $g_{tot}(x)$ increases in a way that reflects the ``{\bf average}" number of Goldbach pairs for a set of integer values of similar magnitude. Therefore the value of its derivative approximates this average value. To obtain the real value we must understand how the number of Goldbach pairs for a given $N$ depends on its factorization into primes, and how this number compares with the average we just mentioned. 
\subsection{\normalsize The $\boldsymbol {N\!DF}$ (divisor factor for $\boldsymbol N$)}
To quantify this dependency, we compare two consecutive $N$ values. We choose the first one ($P$) to be a prime number, while the second ($N_{1} = P+1$) is composite. Suppose the factorization of the $N_{1}$ contains a few odd primes ($p_{1} \dots p_{n}$) smaller than $\sqrt{N}$.
Then, for every of these $p_{i}$ primes, among all the pairs $P-m,P+m$, in 2 out of $p_{i}$ cases either $P-m$ or $P+m$ divide $p_{i}$, leaving a fraction of $(p_{i}-2)/p_{i}$ of possible Goldbach pairs.
Instead, for $N_{1}$, in 1 out of $p_{i}$ cases both $N-m$ and $N+m$ divide $p_{i}$, and the fraction of possible Goldbach pairs which is left is  $(p_{i}-1)/p_{i}$.
So, for each odd prime $p$ in the factorization of a generic N, the proportion of Goldbach pair increases by roughly a factor $(p-1)/(p-2)$ compared with the case if $N$ were not a multiple of $p$.
\\
We call the product over the prime factors of $N$
\begin{equation}
1 \cdot \prod_{\substack {p = 3 \\ N \!\!\!\mod p \,= \, 0}}^{\sqrt{N}} \frac{p_{i}-1}{p_{i}-2}
\end{equation} 
the ``divisor factor for $N$", or {\bf {\em NDF}}.
By definition, the {\em NDF} of a prime number, or of a prime times a power of 2, is 1. If $N$ is a multiple of 3, its {\em NDF} is multiplied by 2 ( i.e. (3-1)/(3-2)). If $N$ is a multiple of 5 (7), the {\em NDF} is multiplied by 4/3 (6/5). It is evident that the smaller a factor of $N$ is, the larger its impact on the {\em NDF} (and on the number of Goldbach pairs) for $N$ is.
\\The max. {\em NDF} value is unbounded, but it grows very slowly with $N$. It increases only when $N$ is a ``{\em primorial}" number, i.e. a product of all consecutive odd primes from 3 onwards. The first few max. {\em NDF} values are $2, \,2.6{\overline 6}, \,3.2, \,3.5{\overline 5}$.
\\More interesting for us is the average $\overline{N\!DF}$ value over a large number of $N$ values. Once we have the value of equation \ref{eq:derivative} for a $N$ value with a known {\em NDF}, we will use $\overline{NDF}$ to determine a scaling factor by which the value has to be multiplied to get an estimate for the number of Goldbach pairs for $N$.
\\To obtain the value of $\overline{N\!DF}$ we computed the {\em NDF} values for a large number of consecutive integers, typically one to ten millions, in different ranges (the highest ones starting at 400000000 and at 9500000000). We then made the average of all these {\em NDF} values.
The result seems to converge to a number close to 1.51477, which is the reciprocal of the Twin Prime constant (0.6601618\dots).\footnote{ Perhaps this result had already been formally proved; but I am unable to find any reference to it.}
\\As we are going to use this result soon, we highlight it:
\begin{equation}
\overline{N\!DF} = \frac{1}{Twin\; Prime\; constant} \approxeq 1.51477
\end{equation}

\subsection{\normalsize Putting all the pieces together}
We are now going to use the different results described above to produce an estimate for the number of Goldbach pairs for $N$ (i.e. pairs of primes {\em N-m, N+m}). Before doing it, let's discuss for a moment how the quality of the approximation \ref{eq:approx} depends on the degree of validity of assumptions a) and b) discussed in \ref{sec-approx}.
\\If we work in the full range (i.e. with $m$ from $0$ to $N-3$), we know that the number of primes from 0 to $N$ will be significantly higher than the number of primes from $N$ to $2N$. Because of this, the number of prime pairs whose sum is not larger than $2N$ will be larger than $P^{2}/4$, and  $G_{tot}(2N)$ will underestimate the correct value. The violation of b) also goes in the same direction, as in the region close to 0 the prime numbers are much denser than in the region close to $2N$. Therefore the proportion of pairs for which one of the primes is very low will be higher than normal, as well as the probability that for such pairs the sum of the two members is lower than $2N$.
\\On the other hand, if we work in the reduced range (and $N$ is large enough), a) and b) are ``almost" respected. There will still be a difference, but in this case we can expect that the estimate will be closer to the real result.
\\
These are the estimates for the number of Goldbach pairs for $N$, respectively in the full and in the reduced range
\begin{equation}
\label{eq:full}
{\rm G\_est\_full}(N) = g_{tot}'(2N) \cdot 2 \cdot 0.66016 \cdot N\!DF(N)
\end{equation}
and
\begin{equation}
\label{eq:reduced}
{\rm G\_est\_reduced}(N) = g_{tot}'(2N) \cdot 2 \cdot 0.66016 \cdot N\!DF(N) \cdot (\sqrt{2}-1)
\end{equation}
The factor 2 before the Twin Prime constant is there because, while the function $G_{tot}$, its approximation $g_{tot}$ and the derivative are defined for all integers, the Goldbach pairs are only distributed among the even integers. The term 0.66016 compensates the fact that the value of $\overline{N\!DF}$, implicitly contained in $G_{tot}$, $g_{tot}$ and $g'_{tot}$, is 1/0.66016.
\newpage
\section{How good is the estimate?}
To check the accuracy of the estimate, we compare it with the exact number of Goldbach pairs for several $N$ values ranging from $5 \cdot 10^{6}$ to $5 \cdot 10^{9}$. For every $N$, we will check both the full and the reduced ranges.
\\
The first table refers to the full range.\\ \ \\
\begin{scriptsize}
\noindent
\makebox[\linewidth] {
\begin{tabular}{|l|l|l|l|l|} \hline
$N$         & $NDF$      & exact count   & estimate        &  estimate/exact         \\  \hline
5000000  &   1.3333  &   38807  &   36317  &   0.9358     \\  \hline
5000001  &   2.0444  &   59624  &   55686  &   0.9339     \\  \hline
5000002  &   1.2706  &   36850  &   34608  &   0.9392     \\  \hline
5000003  &   1.0238  &   29835  &   27886  &   0.9347     \\  \hline
5000004  &   2.0000  &   58229  &   54475  &   0.9355     \\  \hline
5000005  &   1.3468  &   39045  &   36684  &   0.9395     \\  \hline
5000006  &   1.2318  &   35731  &   33550  &   0.9390     \\  \hline
5000007  &   2.0113  &   58445  &   54783  &   0.9373     \\  \hline
5000008  &   1.1024  &   31905  &   30026  &   0.9411     \\  \hline
5000009  &   1.2193  &   35420  &   33210  &   0.9376     \\  \hline
5000010  &   2.6667  &   77536  &   72634  &   0.9368     \\  \hline
5000011  &   1.0000  &   29033  &   27238  &   0.9382     \\  \hline
----     &    ----   &   ----   &   ----   &   ----       \\  \hline
50000000  &   1.3333  &   291400  &   275488  &   0.9454     \\  \hline
50000001  &   2.1223  &   464621  &   438503  &   0.9438     \\  \hline
50000002  &   1.1313  &   247582  &   233747  &   0.9441     \\  \hline
50000003  &   1.0020  &   218966  &   207038  &   0.9455     \\  \hline
50000004  &   2.0000  &   437717  &   413232  &   0.9441     \\  \hline
50000005  &   1.4815  &   323687  &   306098  &   0.9457     \\  \hline
50000006  &   1.2000  &   263241  &   247939  &   0.9419     \\  \hline
50000007  &   2.0000  &   437518  &   413232  &   0.9445     \\  \hline
50000008  &   1.0105  &   220846  &   208791  &   0.9454     \\  \hline
50000009  &   1.0679  &   233634  &   220651  &   0.9444     \\  \hline
50000010  &   2.7259  &   595554  &   563220  &   0.9457     \\  \hline
----     &    ----   &   ----   &   ----   &   ----       \\  \hline
500000000  &   1.3333  &   2274205  &   2161238  &   0.9503     \\  \hline
500000001  &   2.0509  &   3496205  &   3324312  &   0.9508     \\  \hline
500000002  &   1.0256  &   1747858  &   1662491  &   0.9512     \\  \hline
500000003  &   1.0000  &   1704301  &   1620929  &   0.9511     \\  \hline
500000004  &   2.4371  &   4151660  &   3950384  &   0.9515     \\  \hline
500000005  &   1.4222  &   2422662  &   2305321  &   0.9516     \\  \hline
500000006  &   1.1506  &   1960129  &   1864965  &   0.9514     \\  \hline
500000007  &   2.2034  &   3752836  &   3571587  &   0.9517     \\  \hline
500000008  &   1.0000  &   1704555  &   1620929  &   0.9509     \\  \hline
500000009  &   1.0000  &   1703977  &   1620929  &   0.9513     \\  \hline
500000010  &   2.8297  &   4821673  &   4586817  &   0.9513     \\  \hline
----     &    ----   &   ----   &   ----   &   ----       \\  \hline
4900000000  &   1.6000  &   21437787  &   20508113  &   0.9566     \\  \hline
4900000001  &   1.0909  &   14619657  &   13982804  &   0.9564     \\  \hline
4900000002  &   2.0000  &   26804478  &   25635141  &   0.9564     \\  \hline
4900000003  &   1.0593  &   14195586  &   13577612  &   0.9565     \\  \hline
4900000004  &   1.0000  &   13401887  &   12817571  &   0.9564     \\  \hline
4900000005  &   2.6667  &   35754539  &   34180189  &   0.9560     \\  \hline
4900000006  &   1.1111  &   14890471  &   14241745  &   0.9564     \\  \hline
4900000007  &   1.2000  &   16083078  &   15381085  &   0.9564     \\  \hline
4900000008  &   2.0211  &   27088671  &   25904985  &   0.9563     \\  \hline
4900000009  &   1.0169  &   13626809  &   13034818  &   0.9566     \\  \hline
4900000010  &   1.3333  &   17872601  &   17090094  &   0.9562     \\  \hline
\end{tabular}
}
\end{scriptsize}
\\ \ \\
The last column shows by how much the estimate \ref{eq:full} underestimates the exact result. One can notice that 1) the error only depends by the size of $N$ (it does not depend by the $NDF$), and 2) the error decreases when $N$ gets larger.
\newpage
\noindent
In the next table we look at the same values for $N$, but this time we only count the Goldbach pairs in the reduced range (this time the estimate comes from \ref{eq:reduced}).\\  \  \\
\begin{scriptsize}
\noindent
\makebox[\linewidth] {
\begin{tabular}{|l|l|l|l|l|} \hline
$N$         & $NDF$      & exact count   & estimate        &  estimate/exact         \\  \hline
5000000  &   1.3333  &   15378  &   15043  &   0.9782     \\  \hline
5000001  &   2.0444  &   23696  &   23066  &   0.9734     \\  \hline
5000002  &   1.2706  &   14601  &   14335  &   0.9818     \\  \hline
5000003  &   1.0238  &   11881  &   11551  &   0.9722     \\  \hline
5000004  &   2.0000  &   23203  &   22564  &   0.9725     \\  \hline
5000005  &   1.3468  &   15542  &   15195  &   0.9777     \\  \hline
5000006  &   1.2318  &   14176  &   13897  &   0.9803     \\  \hline
5000007  &   2.0113  &   23220  &   22692  &   0.9773     \\  \hline
5000008  &   1.1024  &   12597  &   12437  &   0.9873     \\  \hline
5000009  &   1.2193  &   14145  &   13756  &   0.9725     \\  \hline
5000010  &   2.6667  &   30848  &   30086  &   0.9753     \\  \hline
5000011  &   1.0000  &   11521  &   11282  &   0.9793     \\  \hline
----     &    ----   &   ----   &   ----   &   ----       \\  \hline
50000000  &   1.3333  &   116409  &   114111  &   0.9803     \\  \hline
50000001  &   2.1223  &   185554  &   181634  &   0.9789     \\  \hline
50000002  &   1.1313  &   99058  &   96821  &   0.9774     \\  \hline
50000003  &   1.0020  &   87620  &   85758  &   0.9788     \\  \hline
50000004  &   2.0000  &   174617  &   171166  &   0.9802     \\  \hline
50000005  &   1.4815  &   129281  &   126790  &   0.9807     \\  \hline
50000006  &   1.2000  &   105135  &   102700  &   0.9768     \\  \hline
50000007  &   2.0000  &   174610  &   171166  &   0.9803     \\  \hline
50000008  &   1.0105  &   87999  &   86484  &   0.9828     \\  \hline
50000009  &   1.0679  &   93499  &   91397  &   0.9775     \\  \hline
50000010  &   2.7259  &   237877  &   233293  &   0.9807     \\  \hline
----     &    ----   &   ----   &   ----   &   ----       \\  \hline
500000000  &   1.3333  &   912410  &   895214  &   0.9812     \\  \hline
500000001  &   2.0509  &   1403942  &   1376975  &   0.9808     \\  \hline
500000002  &   1.0256  &   701093  &   688626  &   0.9822     \\  \hline
500000003  &   1.0000  &   683612  &   671411  &   0.9822     \\  \hline
500000004  &   2.4371  &   1666306  &   1636302  &   0.9820     \\  \hline
500000005  &   1.4222  &   972309  &   954895  &   0.9821     \\  \hline
500000006  &   1.1506  &   786276  &   772494  &   0.9825     \\  \hline
500000007  &   2.2034  &   1506137  &   1479400  &   0.9822     \\  \hline
500000008  &   1.0000  &   683901  &   671411  &   0.9817     \\  \hline
500000009  &   1.0000  &   684353  &   671411  &   0.9811     \\  \hline
500000010  &   2.8297  &   1934530  &   1899922  &   0.9821     \\  \hline
----     &    ----   &   ----   &   ----   &   ----       \\  \hline
4900000000  &   1.6000  &   8631392  &   8494739  &   0.9842     \\  \hline
4900000001  &   1.0909  &   5888002  &   5791867  &   0.9837     \\  \hline
4900000002  &   2.0000  &   10795934  &   10618423  &   0.9836     \\  \hline
4900000003  &   1.0593  &   5716500  &   5624031  &   0.9838     \\  \hline
4900000004  &   1.0000  &   5396255  &   5309212  &   0.9839     \\  \hline
4900000005  &   2.6667  &   14403395  &   14157898  &   0.9830     \\  \hline
4900000006  &   1.1111  &   5998489  &   5899124  &   0.9834     \\  \hline
4900000007  &   1.2000  &   6478102  &   6371054  &   0.9835     \\  \hline
4900000008  &   2.0211  &   10907743  &   10730196  &   0.9837     \\  \hline
4900000009  &   1.0169  &   5488978  &   5399198  &   0.9836     \\  \hline
4900000010  &   1.3333  &   7198391  &   7078949  &   0.9834     \\  \hline
\end{tabular}
}
\end{scriptsize}
\\ \ \\
From the table, we can see that the estimate's error is much lower in this case. The reason is that we are examining a region in which the density of primes is much more uniform. In the next section we will look in detail at this problem, and implement a correction to produce an estimate much closer to the reality.
\section{Can we get a better estimate?}
We have seen that, for a given $N$, our formulas constantly underestimate the real number of Goldbach pairs $N-m$, $N+m$, and that the formula for the reduced range is more precise than the one for the full range. The main reason for this is the unbalance between the number of primes below $N$ and the number of primes from $N$ to $2N$. The larger the proportion of primes below $N$ is, the larger the estimate error, independently on the specific value of $N$. This unbalance decreases slowly when $N$ increases. To get a better estimate, we must evaluate how the function that counts the total number of Goldbach pairs depend on the unbalance.
\\
We define the unbalance $U(N)$ as 
\begin{equation}
U(N) = \frac{2 \cdot \pi(N)}{\pi(2N)}
\end{equation}
The value of $U(N)$ is always above 1, and tends to 1 from above when N goes to infinity.
\\In the presence of unbalance (i.e. always), our first approximation (eq. \ref{eq:approx}) underestimates the total number of Goldbach pairs, and the size of the error mainly depends on the value of $U(N)$. We look for a correction in the form $U(N)^{k}$, for some value of $k$, such that if we multiply our first approximation by the correction, we obtain a more precise approximation for the total number of Goldbach pairs. If such a simple correction exists, we can then apply it to the estimate of Goldbach pairs for a given $N$, to get a more precise result.
We start by looking at three values for $k$: 1, 2, and 3/2.
\\The table below reports, for a few values of $N$, the total number of Goldbach pairs up to $2N$, the approximation from \ref{eq:approx}, the ratio between the real value and the approximation, and the values $U(N)$, $U(N)^{2}$ and $U(N)^{3/2}$. 
\\ \ \\
\begin{scriptsize}
\noindent
\makebox[\linewidth] {
\begin{tabular}{|l|l|l|l|l|l|l|} \hline
$N$         & prime pairs total   & estimate        &  total/estimate & $U(N)$ & $U(N)^{2}$ & $U(N)^{3/2}$        \\  \hline
1000000     & 1671879782          & 1540484001      &      {\bf 1.0853}     & 1.0583 & 1.1200     & {\bf 1.0887}  \\  \hline
10000000    & 118268797136        & 110416311810    &      {\bf 1.0711}     & 1.0488 & 1.1000     & {\bf 1.0741}  \\  \hline
100000000   & 8804091976098       & 8298590929256   &      {\bf 1.0609}     & 1.0418 & 1.0853     & {\bf 1.0633}  \\  \hline
1000000000  & 680858394988085     & 646367928470289 &      {\bf 1.0534}     & 1.0367 & 1.0747     & {\bf 1.0555}  \\  \hline
\end{tabular}
}
\end{scriptsize}
\\ \ \\
From the above table, we can say that $U(N)^{3/2}$ provides a satisfactory scaling correction for the unbalance effect (this effect will anyway decrease when N becomes very large, and will tend to 0 for N going to infinity). We will therefore apply this correction to the Goldbach pairs estimate for a given N value.
\\In the next tables (the first for the full range, the second for the reduced one) we take again some of the N values already examined. We add some columns, to show the value of $U(N)$, the correction (i.e. $U(N)^{3/2}$), and the ratio between the corrected estimate and the exact count.\\
In the tables, the columns ``{\em estimate}" and ``{\em estimate/exact}" still refer to the original, uncorrected estimate. The column ``$U(N)^{3/2}$" shows the scaling correction by which the estimate has to be multiplied, and the last column shows the ratio between the corrected estimate and the exact count. As one can see, the corrected estimate is much more accurate. 
\newpage
\begin{scriptsize}
\noindent
\makebox[\linewidth] {
\begin{tabular}{|l|l|l|l|l|l|l|} \hline
$N$         &  exact count   & estimate        &  estimate/exact   & $U(N)$ & $U(N)^{3/2}$ (correction) & {\bf corr. estimate/exact}     \\  \hline
5000000  &   38807  &   36317  &   0.9358  &   1.0488  &  1.0741  &   {\bf 1.0052}    \\  \hline
5000001  &   59624  &   55686  &   0.9339  &   1.0488  &  1.0741  &   {\bf 1.0032}    \\  \hline
5000002  &   36850  &   34608  &   0.9392  &   1.0488  &  1.0741  &   {\bf 1.0088}    \\  \hline
5000003  &   29835  &   27886  &   0.9347  &   1.0488  &  1.0741  &   {\bf 1.0040}    \\  \hline
5000004  &   58229  &   54475  &   0.9355  &   1.0488  &  1.0741  &   {\bf 1.0049}    \\  \hline
5000005  &   39045  &   36684  &   0.9395  &   1.0488  &  1.0741  &   {\bf 1.0092}    \\  \hline
----     &    ----   &  ----   & ----   &  ----   &   ----   &   ----       \\  \hline
50000000  &   291400  &   275488  &   0.9454  &   1.0418  &  1.0633  &   {\bf 1.0053}    \\  \hline
50000001  &   464621  &   438503  &   0.9438  &   1.0418  &  1.0633  &   {\bf 1.0036}    \\  \hline
50000002  &   247582  &   233747  &   0.9441  &   1.0418  &  1.0633  &   {\bf 1.0039}    \\  \hline
50000003  &   218966  &   207038  &   0.9455  &   1.0418  &  1.0633  &   {\bf 1.0054}    \\  \hline
50000004  &   437717  &   413232  &   0.9441  &   1.0418  &  1.0633  &   {\bf 1.0039}    \\  \hline
50000005  &   323687  &   306098  &   0.9457  &   1.0418  &  1.0633  &   {\bf 1.0056}    \\  \hline
----     &    ----   &  ----   & ----   &  ----   &   ----   &   ----       \\  \hline
500000000  &   2274205  &   2161238  &   0.9503  &   1.0367  &  1.0555  &   {\bf 1.0031}    \\  \hline
500000001  &   3496205  &   3324312  &   0.9508  &   1.0367  &  1.0555  &   {\bf 1.0036}    \\  \hline
500000002  &   1747858  &   1662491  &   0.9512  &   1.0367  &  1.0555  &   {\bf 1.0039}    \\  \hline
500000003  &   1704301  &   1620929  &   0.9511  &   1.0367  &  1.0555  &   {\bf 1.0039}    \\  \hline
500000004  &   4151660  &   3950384  &   0.9515  &   1.0367  &  1.0555  &   {\bf 1.0043}    \\  \hline
500000005  &   2422662  &   2305321  &   0.9516  &   1.0367  &  1.0555  &   {\bf 1.0044}    \\  \hline
----     &    ----   &  ----   & ----   &  ----   &   ----   &   ----       \\  \hline
4900000000  &   21437787  &   20508113  &   0.9566  &   1.0327  &  1.0494  &   {\bf 1.0039}    \\  \hline
4900000001  &   14619657  &   13982804  &   0.9564  &   1.0327  &  1.0494  &   {\bf 1.0037}    \\  \hline
4900000002  &   26804478  &   25635141  &   0.9564  &   1.0327  &  1.0494  &   {\bf 1.0036}    \\  \hline
4900000003  &   14195586  &   13577612  &   0.9565  &   1.0327  &  1.0494  &  {\bf  1.0037}    \\  \hline
4900000004  &   13401887  &   12817571  &   0.9564  &   1.0327  &  1.0494  &   {\bf 1.0037}    \\  \hline
4900000005  &   35754539  &   34180189  &   0.9560  &   1.0327  &  1.0494  &   {\bf 1.0032}    \\  \hline
\end{tabular}
}
\end{scriptsize}
\\Results for the full range
\\ \ \\
\begin{scriptsize}
\noindent
\makebox[\linewidth] {
\begin{tabular}{|l|l|l|l|l|l|l|} \hline
$N$         &  exact count   & estimate        &  estimate/exact   & $U(N)$ & $U(n)^{3/2}$ (correction) & {\bf corr. estimate/exact}     \\  \hline
5000000  &   15378  &   15043  &   0.9782  &   1.0142  &  1.0213  &  {\bf  0.9991}    \\  \hline
5000001  &   23696  &   23066  &   0.9734  &   1.0142  &  1.0213  &  {\bf  0.9942}    \\  \hline
5000002  &   14601  &   14335  &   0.9818  &   1.0142  &  1.0213  &  {\bf  1.0027}    \\  \hline
5000003  &   11881  &   11551  &   0.9722  &   1.0142  &  1.0213  &  {\bf  0.9929}    \\  \hline
5000004  &   23203  &   22564  &   0.9725  &   1.0142  &  1.0213  &  {\bf  0.9932}    \\  \hline
5000005  &   15542  &   15195  &   0.9777  &   1.0142  &  1.0213  &  {\bf  0.9985}    \\  \hline
----     &    ----   &  ----   & ----   &  ----   &   ----   &   ----       \\  \hline
50000000  &   116409  &   114111  &   0.9803  &   1.0121  &  1.0183  &  {\bf  0.9982}    \\  \hline
50000001  &   185554  &   181634  &   0.9789  &   1.0121  &  1.0183  &  {\bf  0.9967}    \\  \hline
50000002  &   99058  &   96821  &   0.9774  &   1.0121  &  1.0183  &   {\bf 0.9953}    \\  \hline
50000003  &   87620  &   85758  &   0.9788  &   1.0121  &  1.0183  &   {\bf 0.9966}    \\  \hline
50000004  &   174617  &   171166  &   0.9802  &   1.0121  &  1.0183  &  {\bf  0.9981}    \\  \hline
50000005  &   129281  &   126790  &   0.9807  &   1.0121  &  1.0183  &  {\bf  0.9986}    \\  \hline
----     &    ----   &  ----   & ----   &  ----   &   ----   &   ----       \\  \hline
500000000  &   912410  &   895214  &   0.9812  &   1.0107  &  1.0161  &   {\bf 0.9970}    \\  \hline
500000001  &   1403942  &   1376975  &   0.9808  &   1.0107  &  1.0161  &  {\bf  0.9966}    \\  \hline
500000002  &   701093  &   688626  &   0.9822  &   1.0107  &  1.0161  &   {\bf 0.9980}    \\  \hline
500000003  &   683612  &   671411  &   0.9822  &   1.0107  &  1.0161  &  {\bf  0.9980}    \\  \hline
500000004  &   1666306  &   1636302  &   0.9820  &   1.0107  &  1.0161  &  {\bf  0.9978}    \\  \hline
500000005  &   972309  &   954895  &   0.9821  &   1.0107  &  1.0161  &   {\bf 0.9979}    \\  \hline
----     &    ----   &  ----   & ----   &  ----   &   ----   &   ----       \\  \hline
4900000000  &   8631392  &   8494739  &   0.9842  &   1.0096  &  1.0144  &  {\bf  0.9984}    \\  \hline
4900000001  &   5888002  &   5791867  &   0.9837  &   1.0096  &  1.0144  &  {\bf  0.9979}    \\  \hline
4900000002  &   10795934  &   10618423  &   0.9836  &   1.0096  &  1.0144  &  {\bf  0.9977}    \\  \hline
4900000003  &   5716500  &   5624031  &   0.9838  &   1.0096  &  1.0144  &   {\bf 0.9980}    \\  \hline
4900000004  &   5396255  &   5309212  &   0.9839  &   1.0096  &  1.0144  &   {\bf 0.9981}    \\  \hline
4900000005  &   14403395  &   14157898  &   0.9830  &   1.0096  &  1.0144  &  {\bf  0.9971}    \\  \hline
\end{tabular}
}
\end{scriptsize}
\\Results for the reduced range
\\
\section{Conclusions}
Using only simple concepts, and without relying on the exact distribution of the prime numbers, we have derived an estimate for the number of Goldbach pairs {\em N-m, N+m} adding to $2N$. A first, rough version of the estimate only requires the knowledge of the number of primes up to $2N$, of the average value of the {\em NDF} (the divisor factor for $N$), and of the prime factorization of $N$. The knowledge of the ratio between the number of primes from 3 to $N$ and the number of primes from $N$ to $2N-3$ allows us to determine a ``scaling correction". When the first estimate is multiplied by the scaling correction, we obtain a very accurate estimate.

\bibliography{Primes}
\bibliographystyle{unsrt} 
\end{document}